\documentclass{article}
\usepackage{amsfonts,amsmath,amsthm,amssymb,graphics,epsfig}

\begin{document}

\title{SINGULAR POINTS OF REAL QUARTIC CURVES VIA COMPUTER ALGEBRA}
\bigskip
\author{David A. Weinberg and Nicholas J. Willis}

\date{}
\maketitle
\begin{abstract}
There are thirteen types of singular points for irreducible real
quartic curves and seventeen types of singular points for
reducible real quartic curves.  This classification is originally
due to D.A. Gudkov. There are nine types of singular points for
irreducible complex quartic curves and ten types of singular
points for reducible complex quartic curves. We derive the
complete classification with proof by using the computer algebra
system Maple.  We clarify that the classification is based on
computing just enough of the Puiseux expansion to separate the
branches. Thus, the proof consists of a sequence of large symbolic
computations that can be done nicely using Maple.
\end{abstract}

\section{Introduction}
The classification of singular points of real quartic curves is originally due to D.A. Gudkov [2,3,4,5].
He determined the individual types of singular points, as well as all possible sets of singular points that real
quartic curves can have.  In this paper, we will derive the thirteen individual types of singular points for
irreducible real quartic curves and the seventeen individual types of singular points for reducible real quartic curves.
 Though our results are not new, the description of the equivalence relation is new and our proof is new and
  gives a very nice illustration of the role that computer algebra can play in doing
  proofs. Furthermore, our proof is self-contained and is the most
  elementary proof possible, which makes this material accessible
  to the widest possible audience.

The general question is how shall we classify singular points of
real quartic curves.  For each fixed degree $n$, we want a finite
classification of singular points for all algebraic curves of
degree n.  Thus, in general, the local diffeomorphism type is not
the desired criterion of classification for singular points.  For
irreducible quartic curves, there are only finitely many
diffeomorphism types, but for reducible quartic curves, there are
infinitely many equivalence classes with respect to local
diffeomorphism.  For example, in the Arnol'd notation, four lines
intersecting at the origin represents an $X_9$ singular point
which is really an infinite family of smoothly inequivalent
singularities. Notice here that an irreducible real quintic curve
can have an $X_9$ singular point. The tradition is to treat these
as one class by fiat. In our scheme the $X_9$ family will appear
naturally as a single class.

Now let us describe how we will classify the individual types of
singular points that a real quartic curve can have.  Given any
polynomial equation $F(x,y)=0$, it is possible to solve for $y$ in
terms of $x$ in the form of fractional power series, called
Puiseux expansions.  There is an algorithm for doing this, and the
software Maple computes such Puiseux expansions, even for curves
with literal coefficients.  Our classification is based on taking
just enough of the Puiseux expansions to separate the ``branches,"
and noting the exponents at which the ``branches" separate.  In
other words, compute the Puiseux expansions to a power of $x$ such
that all expansions are unique.  Then we will associate a
tree-type graph, to which we will refer as a ``tree diagram" or
``diagram."  These diagrams will be described in detail below and
will codify how the ``branches" separate and will serve to
classify the type of the singular point.  At this point let us
remark that the term ``branch" already has a traditional meaning
in this context.  We are really interested in the distinct Puiseux
expansions.  In [7], C.T.C. Wall has coined the term ``pro-branch"
for the distinct Puiseux expansions.

In studying a singular point of an algebraic curve, the first
thing to look at is the Newton polygon. (Our Newton polygons will
follow the style of Walker [6].)  Corresponding to each segment of
the Newton polygon, there is a quasihomogeneous polynomial [p.
195,1]. If all such quasihomogeneous polynomials have no multiple
factors, then the Newton polygon already tells us the type of the
singularity.  (Note that in this case, we know right away the
exponents at which all of the Puiseux expansions separate.)  But
if there is a multiple factor, then it is necessary to examine the
situation more closely.  For this, we turn to the Puiseux
expansions.  As indicated above, the relevant definition on which
our classification is based is new and appeals to the Puiseux
expansions in an invariant way.

Let us note that we will classify the \emph{real} singular points.
(It is possible for a real quartic curve to have a complex
conjugate pair of singular points.  We will avoid this case.)  By
a simple translation of axes, we may assume that the singular
point is at the origin.  We will treat both irreducible and
reducible curves, but note that the notions of irreducible and
reducible are with respect to the complex numbers.  Note also that
we will not study reducible curves with multiple components.

The objects being classified are pairs whose first coordinate is a
real quartic curve, specified by a 4th degree polynomial with real
coefficients, considered up to a real nonzero multiplicative
constant, and the second coordinate is a singular point of the
curve in the first coordinate.  Let the quartic curve be given by
$f(x,y)=0$, where
\newpage
\begin{align*}
f(x,y)&=a_{00} + a_{10} x + a_{01} y + a_{20} x^2 + a_{11} xy \\
&+a_{02} y^2 + a_{30} x^3 + a_{21} x^2y + a_{12} xy^2 + a_{03} y^3+ a_{40} x^4 + a_{31}x^3 y \\
 &+ a_{22} x^2 y^2 + a_{13} x y^3 + a_{04} y^4.
 \end{align*}

Since we may assume that our singular point is at the origin, we
have $a_{00}=0$.  Since the point is singular, $ a_{10} = a_{01} =
0$.  In this paper we will use the term ``tangent cone" to refer
to the terms of lowest degree in $f(x,y)$. The degree of these
terms is called the multiplicity of the point.  If the point is of
multiplicity four, then the curve must be reducible since any
homogeneous polynomial of degree $4$ must factor.  Thus, for
irreducible curves, we only need to study points of multiplicity
three or two.

Let us now explain how all of the cases are enumerated. First we
choose the tangent cone by choosing the tangent lines together
with their multiplicities.  The choice of tangent lines can be
fixed by a linear change of coordinates.  Moreover, by rotation of
axes, we may assume no tangent line is vertical. For each tangent
cone, we consider all possible Newton polygons. For each Newton
polygon, we first consider the case where none of the
quasihomogeneous polynomials corresponding to the segments of the
Newton polygon have a multiple factor. Then we consider the cases
where there is a multiple factor.  When there is a multiple
factor, the choice of this factor can be fixed by a linear change
of coordinates. For irreducible quartic curves, the only case of
this kind is the one with tangent cone $(y + x^2)^2$.  In this
case, Maple is used to compute the Puiseux expansions for the
corresponding family of curves. The different types of singular
points are then determined by the vanishing or nonvanishing of
certain polynomials in the coefficients of this family; these
polynomials are given to us by the Maple computations in the form
of discriminant-like polynomial coefficients of the Puiseux
expansions. The details of this are carried out in the next
section.

To be more specific, the details of the following outline will be
carried out in the next section.  For irreducible quartic curves,
by a linear change of coordinates, as described above, it suffices
to consider the following families:

Multiplicity 3

$$
y^3 + ax^4 + bx^3y + cx^2y^2 + dxy^3 + ey^4 = 0, \qquad a\neq 0
 $$

 $$
 y^2(y - x)+ ax^4 + bx^3y + cx^2y^2 + dxy^3 + ey^4 = 0, \qquad a\neq 0
 $$

 $$y(y - x)(y - 2x) + ax^4 + bx^3y + cx^2y^2 + dxy^3 + ey^4 = 0, \qquad a\neq 0
 $$

 $$y(y^2 + x^2) + ax^4 + bx^3y + cx^2y^2 + dxy^3 + ey^4 = 0, \qquad a\neq 0
 $$
\\
 Multiplicity 2

 $$y^2 - x^2 + ax^3 + bx^2y + cxy^2 + dy^3 + ex^4 + fx^3y +
 gx^2y^2 + hxy^3 + jy^4 = 0
 $$

 $$y^2 + x^2 + ax^3 + bx^2y + cxy^2 + dy^3 + ex^4 + fx^3y +
 gx^2y^2 + hxy^3 + jy^4 = 0
 $$

 $$y^2 + + ax^3 + bx^2y + cxy^2 + dy^3 + ex^4 + fx^3y +
 gx^2y^2 + hxy^3 + jy^4 = 0, \qquad a\neq 0
 $$

 $$(y + x^2)(y - x^2) + ax^3y + bxy^2 + cx^2y^2 + dy^3 + exy^3 +
 fy^4 = 0
 $$

 $$y^2 + x^4 + ax^3y + bxy^2 + cx^2y^2 + dy^3 + exy^3 +
 fy^4 = 0
 $$

 $$(y + x^2)^2 + ax^3y + bxy^2 + cx^2y^2 + dy^3 + exy^3 +
 fy^4 = 0
 $$

 Maple computation is needed only for the last family above.  An
 interesting feature of the proof occurs at the end of this
 computation, where we show that every curve in the family

 $$(y + x^2)^2 + bx^3y + bxy^2 + (1/4b^2 + d)x^2y^2 + dy^3 +
 1/2bdxy^3 + fy^4 = 0
 $$
is reducible.  This is the key to establishing that the list of
 double points is complete.

 For reducible quartic curves, Maple computation is used to
 examine the cases where an irreducible cubic is tangent to the
 line component and where two conics share a common tangent.  (The
 other cases are enumerated by mathematical common sense.)

Given an algebraic curve with a singular point at the origin, let
us now describe how to associate a tree diagram to this singular
point once we have the Puiseux expansions.  Each time at least one
``branch" separates, record the exponent where that happens. Place
all such exponents in a row at the top.
 For each exponent in the top row, there corresponds a column of vertices.
 Each Puiseux expansion corresponds to exactly one vertex in that column, and those expansions
 with the same coefficients up to that exponent correspond to the same vertex.
 Braces will join those pairs of vertices, within a given column, that correspond to complex conjugate coefficients.
  In such a case, the only real solution of the original equation, satisfying the pair of expansions
  indicated by the braces, in a small enough neighborhood of the origin is (0,0). \\

  In [7] C.T.C. Wall uses the term ``pro-branches" to refer to the
  distinct Puiseux expansions belonging to a given singular point,
  and then defines a notion of \emph{exponent of contact} between
  two pro-branches.  It follows from Lemma 4.1.1 on page 68 of
  [7], that the diagram we assign to a singular point is invariant
  under a linear change of coordinates.

 \underline{Example.}  $y^2 = -x^3$. Notice that $y = \pm\; i\; x^{3/2}$, which can also be written as $y = \pm (-x)^{3/2}$.  For each $x < 0$, there are
two distinct real solutions for $y$.  Hence, the diagram is as
shown below (without braces!).
\begin{figure}[!h]\scalebox{.33}{\includegraphics{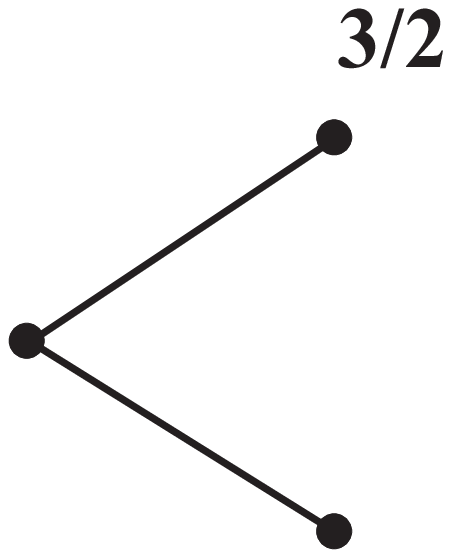}}\end{figure}

We start with one vertex on the left corresponding to the power
zero. Line segments are drawn connecting the vertices from left to
right, where each polygonal path from left to right corresponds to
Puiseux expansions having the same set of coefficients up to a
given exponent. The diagram stops at the first exponent where each
vertex in that column corresponds to exactly one Puiseux
expansion.  The key point is that this tree diagram uniquely
specifies the singularity type (up to permutations of vertices
within columns) provided that no tangent
line at the origin is vertical.\\
\\ \underline{Example.}  $x^2y + x^4 + 2xy^2 + y^3 = 0$.  If $B:= x^2 y + x^4 + 2xy^2 + y^3$, the Maple command puiseux $(B, x = 0,y,3)$ tells us that the Puiseux expansions begin as follows:
$$y = -x+x^{3/2} \qquad (\text{branch} \#1)$$
$$y = -x - x^{3/2} \qquad (\text{branch} \#2)$$
$$y = -x^2 \qquad \qquad \;\;(\text{branch} \#3)$$
In the next section, we will refer to the relevant truncated portion of the Puiseux expansion as the \emph{Puiseux jet}.  Notice that the coefficient of $x$ in branch \#1 and branch \#2 is $-1$, while the coefficient of $x$ in branch \#3 is $0$.  So there is a splitting at the first power of $x$, which is indicated as \\ \begin{figure}[!h]\scalebox{.33}{\includegraphics{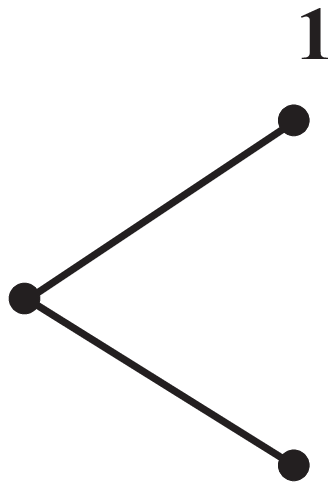}}\end{figure}\\  Next we must show the splitting of \#1 from \#2.  Notice that the power of $x$ at which \#1 and \#2 split is 3/2.
 \newpage Now our diagram looks like the following:\\ \begin{figure}[!h]\scalebox{.33}{\includegraphics{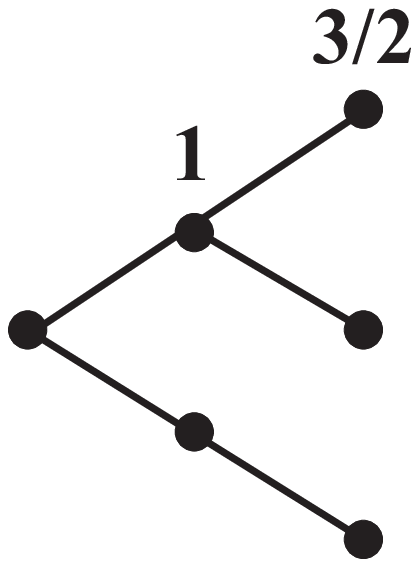}}\end{figure} \\

The diagram is now complete; notice that there are three distinct vertices in the column labeled $3/2$.  \\

\textbf{Definition of the equivalence relation:} two singular
points are equivalent if they have the same diagram as described
above.\\

To summarize, we have described a precise procedure for assigning
a diagram to a singular point of an algebraic curve and this
assignment is invariant under a linear change of coordinates.
\\
\newline
 \underline{Acknowledgment}.  The authors wish to thank Tomas
Recio (Universidad de Cantabria, Santander, Spain), Carlos
Andradas (Universidad Complutense de Madrid, Spain), Eugenii
Shustin (Tel-Aviv University), Jeffrey M. Lee (Texas Tech
University) and Anatoly Korchagin (Texas Tech University) for
several useful conversations.
\\
\newline
 \textbf{Note}: Due to diagram typesetting considerations,
section 2 begins on the following page.
\newpage
\section{Classification and proof.}
\underline{Irreducible curves}.\\ \\
\underline{Multiplicity 3}.\\ \\
Tangent cone:  $y^3$.\\ \\
\begin{table}[th]
\begin{tabular}{l l}
\begin{tabular}{l}Newton polygon: \\ \\ \\ \\  \end{tabular} & \scalebox{.33}{\includegraphics{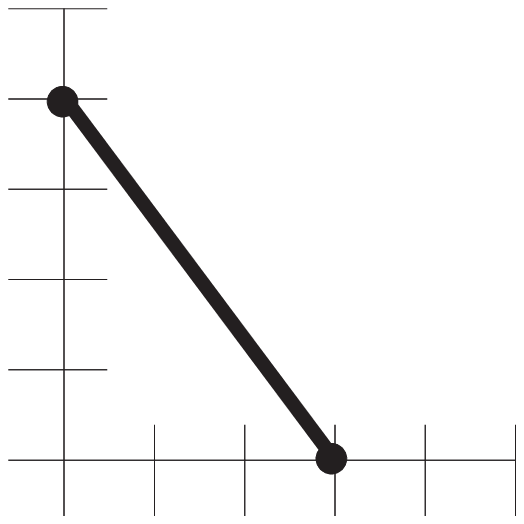}} \end{tabular} \\
$y^3 + ax^4 + bx^3 y + cx^2y^2 + dxy^3 + ey^4 = 0, a \neq 0$.\\ \\
Puiseux jets : $y = -a^{1/3} x^{4/3}$ (three expansions here).\\ \\
\begin{tabular}{l l}
\begin{tabular}{l}
Diagram Type 1:  \\ \\ \\ \\  \end{tabular}&\scalebox{.33}{\includegraphics{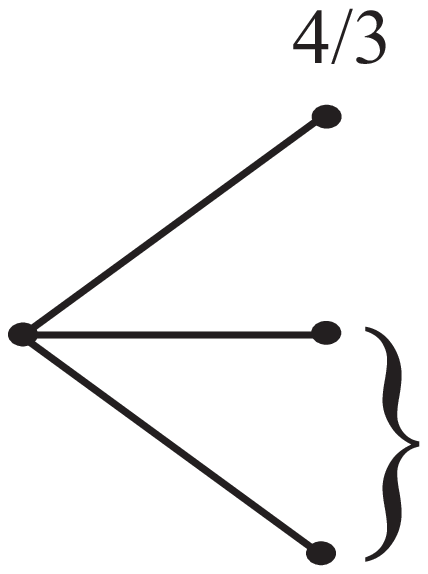}}\end{tabular}\\
\\
Tangent cone:  $y^2 (y-x)$.\\ \\
\begin{tabular}{l l}
\begin{tabular}{l}
Newton polygon: \\ \\ \\ \\ \end{tabular} & \scalebox{.33}{\includegraphics{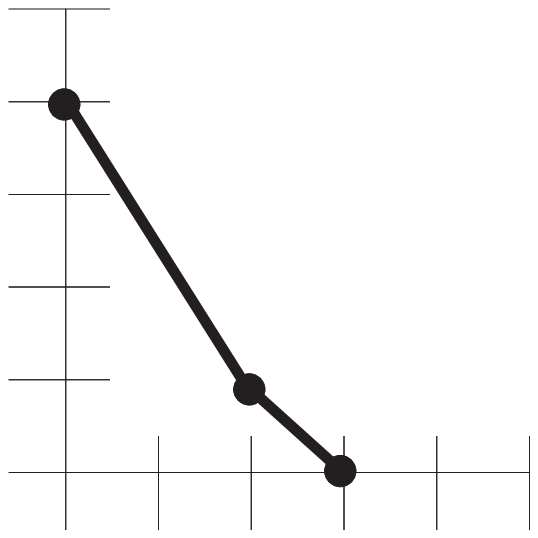}} \end{tabular} \\
$y^2(y-x) + ax^4 + bx^3 y + cx^2y^2 + dxy^3 + ey^4, a \neq 0$.\\ \\
Puiseux jets (from Newton polygon; Maple not needed):\\  $y = x \qquad  \\y = \pm \sqrt{a} x^{3/2}$.\\ \\ \end{table}
\newpage
\begin{table}[t]
\begin{tabular}{l l}
\begin{tabular}{l}Diagram Type 2: \\ \\ \\ \\ \end{tabular}&\scalebox{.33}{\includegraphics{4fig3c}} \end{tabular} \\
\\
Tangent cone:  $y(y - x) (y - 2x)$.\\ \\
\begin{tabular}{l l}
\begin{tabular}{l}Newton polygon: \\ \\ \\ \\ \end{tabular}&\scalebox{.33}{\includegraphics{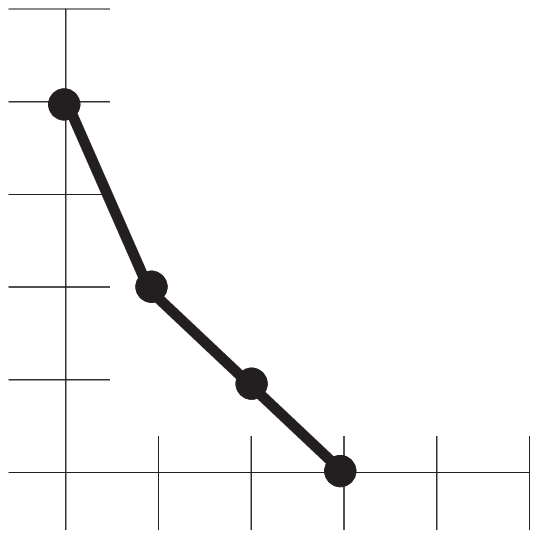}} \end{tabular}\\ \\
$y (y - x) (y - 2x) + ax^4 + bx^3 y + cx^2 y^2 + d x y^3 + e y^4 = 0, a \neq 0$.\\ \\
Puiseux jets:\\
$y = 0$\\
$y = x$\\
$y = 2x$.\\
\begin{tabular}{l l}
\begin{tabular}{l}Diagram Type 3: \\ \\ \\ \\ \end{tabular}& \scalebox{.33}{\includegraphics{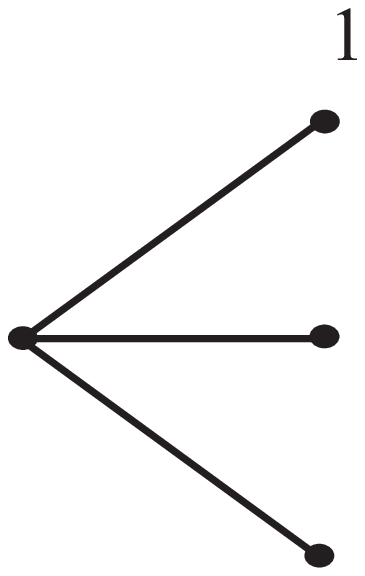}} \end{tabular} \\
\\
Tangent cone: $y (y^2+x^2)$.\\ \\
\begin{tabular}{l l}
\begin{tabular}{l}
Newton polygon: \\ \\ \\ \\ \end{tabular}& \scalebox{.33}{\includegraphics{4fig8c}}\end{tabular} \\ \\
$y(y^2+x^2) + ax^4 + bx^3y + cx^2y^2 + dx y^3 + e y^4 = 0, a \neq 0.$\\ \\
Puiseux jets:\\
$y = 0$\\
$y = \pm i x$.\\ \end{table}
\newpage
\begin{table}[t]
\begin{tabular}{l l}
\begin{tabular}{l}
Diagram Type 4: \\ \\ \\ \\ \end{tabular} & \scalebox{.33}{\includegraphics{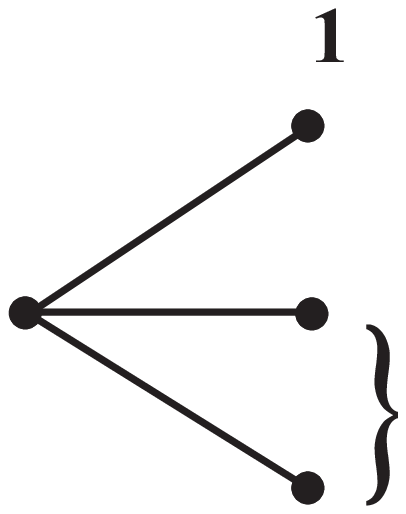}} \end{tabular}\\ \\
\underline{Multiplicity 2.}\\ \\
Tangent cone:  $(y - x) (y + x) = y^2 - x^2$.\\ \\
\begin{tabular}{l l}
\begin{tabular}{l}
Newton polygon: \\ \\ \\ \\ \end{tabular}& \scalebox{.33}{\includegraphics{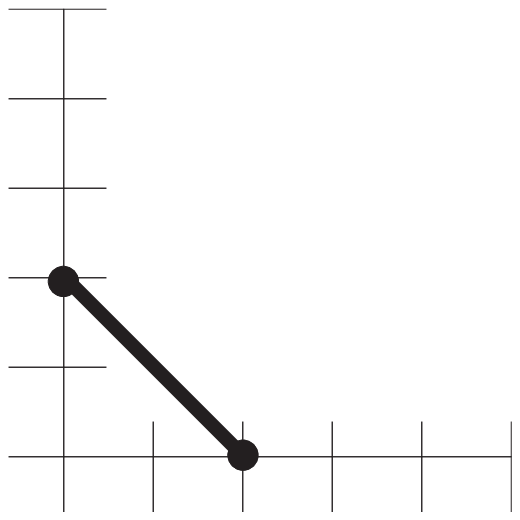}}\end{tabular}\\ \\
$y^2 - x^2 + ax^3 + b x^2 y + cxy^2 + dy^3 + e x^4 + fx^3 y + g x^2 y^2 + hxy^3 + jy^4 = 0$.\\ \\
Puiseux jets:\\
$y = x$\\
$ y= -x$.\\ \\
\begin{tabular}{l l}
\begin{tabular}{l}
Diagram Type 5: \\ \\ \\ \\ \end{tabular} &\scalebox{.33}{\includegraphics{4fig2c}}\end{tabular}\\ \\
Tangent cone: $y^2 + x^2$.\\ \\
\begin{tabular}{l l}
\begin{tabular}{l}
Newton polygon: \\ \\ \\ \\ \end{tabular}& \scalebox{.33}{\includegraphics{4fig12c}} \end{tabular}\\ \\
$y^2+x^2+ax^3+bx^2y +cxy^2 + dy^3 + ex^4 + fx^3y + gx^2y^2+ hxy^3 + jy^4 = 0$.\\ \\
Puiseux jets:\\
$y = \pm i x$.\\ \\
 \end{table}
\newpage
\begin{table}[t]
\begin{tabular}{l l}
\begin{tabular}{l}
Diagram Type 6:  \\ \\ \\ \\ \end{tabular}&\scalebox{.33}{\includegraphics{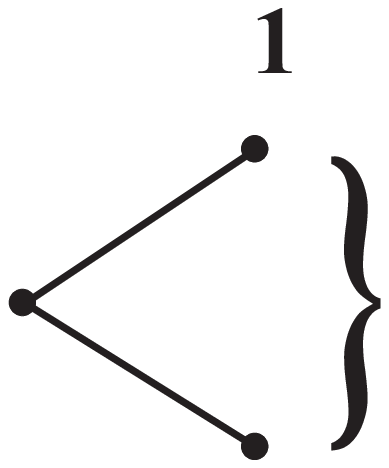}} \end{tabular}\\ \\
Tangent cone: $y^2$.\\ \\
\begin{tabular}{l l}
\begin{tabular}{l}
Newton polygon: \\ \\ \\ \\ \end{tabular} &\scalebox{.33}{\includegraphics{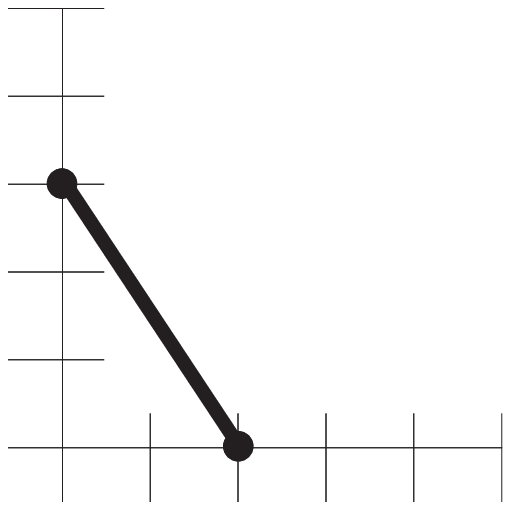}}\end{tabular}\\ \\
$y^2 + ax^3 + bx^2y + cxy^2 + dy^3 + ex^4 + fx^3y + gx^2y^2+hxy^3 + jy^4 = 0, a \neq 0$.\\ \\
Puiseux jets:\\
$y = \pm \sqrt{a}\; x^{3/2}$.\\ \\
\begin{tabular}{l l}
\begin{tabular}{l}
Diagram Type 7: \\ \\ \\ \\ \end{tabular}&\scalebox{.33}{\includegraphics{4fig1c}}\end{tabular}\\ \\
\begin{tabular}{l l}
\begin{tabular}{l}
Newton polygon: \\ \\ \\ \\ \end{tabular} &\scalebox{.33}{\includegraphics{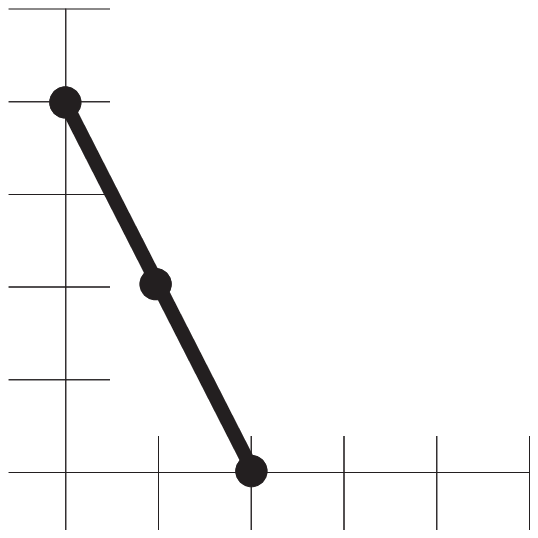}}\end{tabular}\\ \\
Quasihomogeneous factors: $(y + x^2) (y-x^2)$.\\ \\
$(y + x^2)(y - x^2) + ax^3y+bxy^2 + cx^2 y^2 + dy^3 + exy^3 + fy^4 = 0$.\\ \\
Puiseux jets:\\
$y = x^2$\quad \\ $y = -x^2$.\\ \\
\begin{tabular}{l l}
\begin{tabular}{l}
Diagram Type 8: \\ \\ \\ \\
\end{tabular}&\scalebox{.33}{\includegraphics{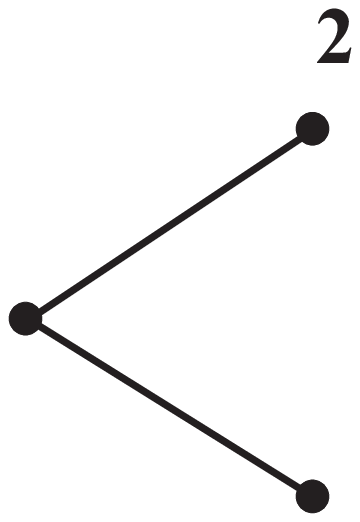}}\end{tabular}\\
\end{table}
\newpage
\begin{table}[t]
Quasihomogeneous factors:  $y^2+x^4$.\\ \\

$y^2+x^4 + ax^3y + bxy^2 + cx^2y^2 + dy^3 + exy^3 + fy^4 = 0$.\\ \\
Puiseux jets:\\
$y = \pm i x^2.$\\ \\
\begin{tabular}{l l}
\begin{tabular}{l}
Diagram Type 9: \\ \\ \\ \\ \end{tabular}&\scalebox{.33}{\includegraphics{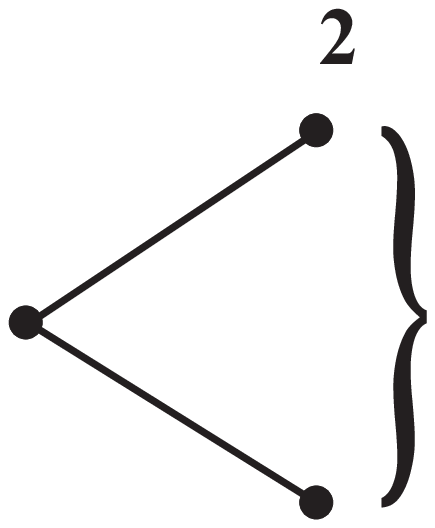}}\end{tabular} \\ \\
Quasihomogeneous factors: $(y + x^2)^2$.\\ \\
$A: = (y + x^2)^2 + ax^3y + bxy^2 + cx^2y^2 + dy^3 + exy^3 + fy^4 = 0$.\\ \\
Notice that the quasihomogeneous polynomial $(y + x^2)^2$ has a double root.  Thus, the family above contains several different types of singular points.  We will determine polynomial conditions on the coefficients that will give all the different types of singular points by using Maple to compute a succession of Puiseux expansions.  We begin by computing the Puiseux expansion of $A$ using the Maple command puiseux$(A, x = 0, y, 0)$.  Notice that the zero in the last argument instructs Maple to exhibit just enough of the Puiseux expansion to separate the ``branches"!\\ \\
Puiseux jets:  $y = -x^2 + (a-b)^{1/2} x^{5/2}$.\\ \\
Condition:  $a \neq b$.\\ \\
\begin{tabular}{l l}
\begin{tabular}{l}
Diagram Type 10: \\ \\ \\ \\ \end{tabular}&\scalebox{.33}{\includegraphics{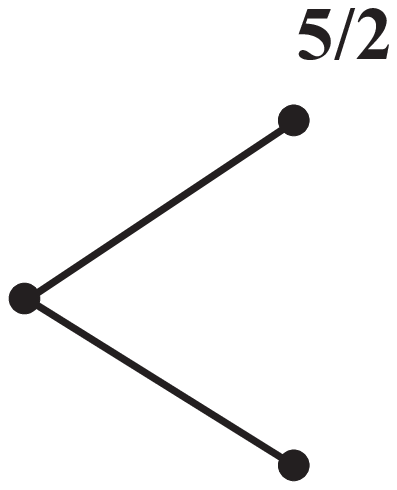}}\end{tabular}\\ \\
Case: $a = b$.\\ \\
Puiseux jets: \\$y = -x^2+x^3$ Root Of $(-b Z + Z^2 + c-d)$.\\ \\
Condition:  $c \neq \frac{1}{4} b^2 + d$.\\ \end{table} \newpage
\begin{table}[t]
\begin{tabular}{l l}
\\ \\
Diagram Type 11. if $b^2 - 4 (c-d) > 0$&  \ Diagram Type 12. if $b^2 - 4 (c-d) < 0$.\\
\scalebox{.33}{\includegraphics{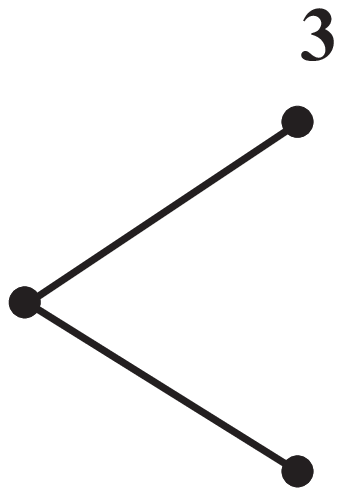}}&\scalebox{.33}{\includegraphics{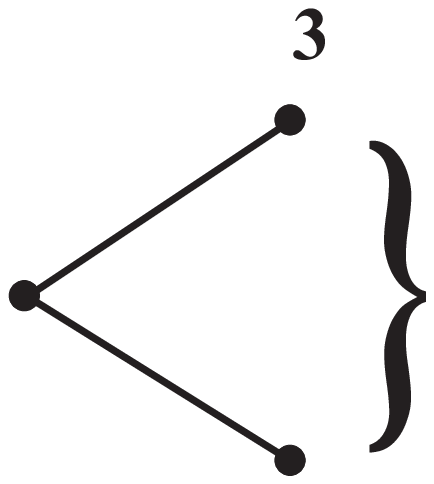}}
\end{tabular} \\ \\
Case: $c = \frac{1}{4} b^2 + d$.\\ \\
Puiseux jets:\\
$y = -x^2 + (e - \frac{1}{2} bd)^{1/2} x^{7/2}$.\\ \\
Condition:  $e \neq \frac{1}{2} bd.$\\ \\
\begin{tabular}{l l}
\begin{tabular}{l}
Diagram Type 13: \\ \\ \\ \\ \end{tabular} &\scalebox{.33}{\includegraphics{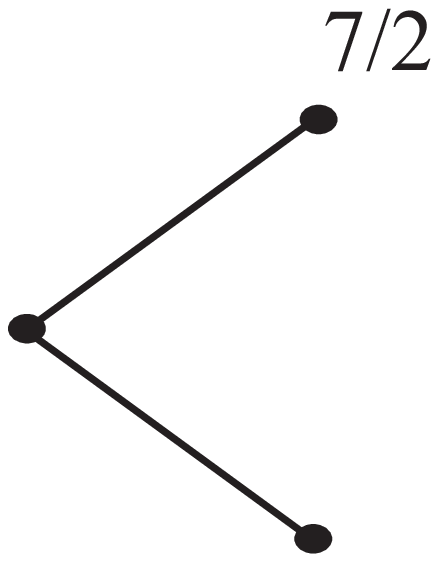}}\end{tabular}\\ \\
Case:  $ e = \frac{1}{2} bd$.\\ \\
Notice that now the family has become\\ \\
$H:= (y + x^2)^2 + bx^3 y + bxy^2 + (\frac{1}{4} b^2 + d)x^2y^2 + dy^3 + \frac{1}{2} bdxy^3 + fy^4$. \\ \\
We cannot show that $H$ is reducible by using the Maple command factor$(H)$. However, the Maple command factor$(H - fy^4)$ shows that $H - fy^4 = \frac{1}{4} (2x^2 + 2y + bxy)(bxy + 2x^2 + 2dy^2 + 2y).$  Therefore, $  H = \frac{1}{4} (2x^2 + bxy + 2y)^2 + \frac{d}{2} (2x^2 + bxy + 2y) y^2 + f y^4$, which is homogeneous in $(2x^2 + bxy + 2y)$ and $y^2$, and thus factors.
\\ \\
This completes the classification of singular point types for irreducible real quartic curves.\\ \end{table} \newpage
\begin{table}[t]
\underline{Reducible Curves}\\ \\
Degrees of factors:  $3,1$.\\ \\
If the straight line does not pass through $(0,0)$, then there are three cases:\\ \\
\begin{tabular}{l l l l}
 && Diagram Type 1: \\
&\scalebox{.33}{\includegraphics{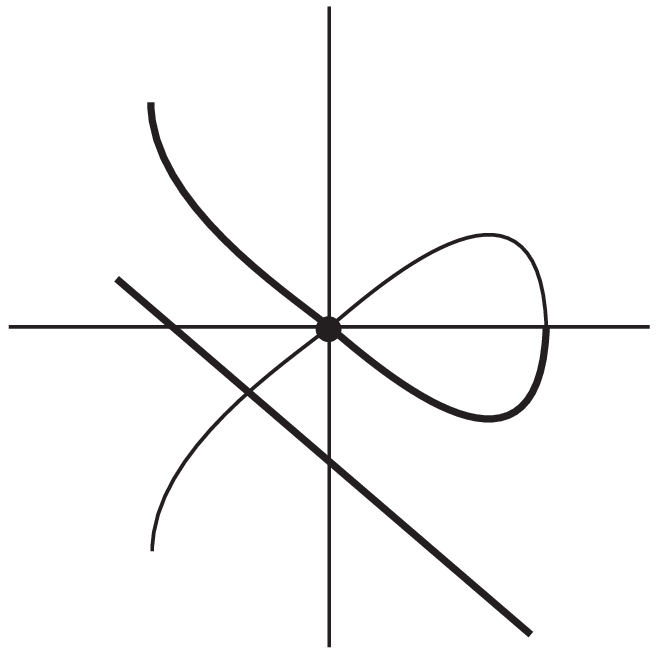}}&&\scalebox{.33}{\includegraphics{4fig2c}} \\ \\
 && Diagram Type 2: \\
 &\scalebox{.33}{\includegraphics{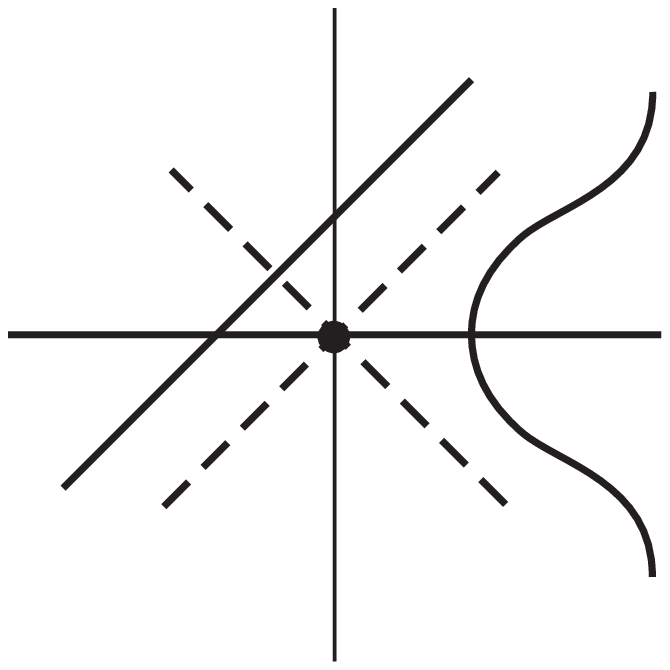}}&&\scalebox{.33}{\includegraphics{4fig15c}}\\ \\
 && Diagram Type 3: \\ &\scalebox{.33}{\includegraphics{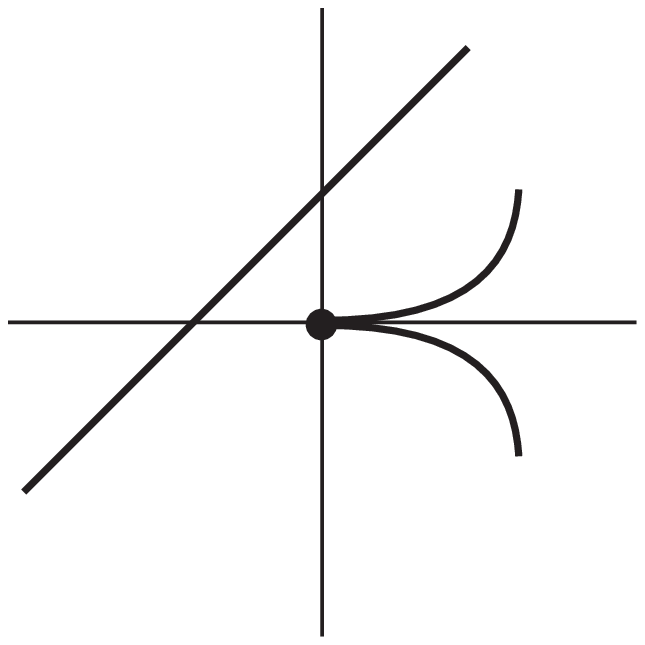}}&&\scalebox{.33}{\includegraphics{4fig1c}}\end{tabular} \\ \\
If the straight line does pass through $(0,0)$, then there are five cases:\\ \\
\begin{tabular}{l l l l}
 && Diagram Type 4: \\ &\scalebox{.33}{\includegraphics{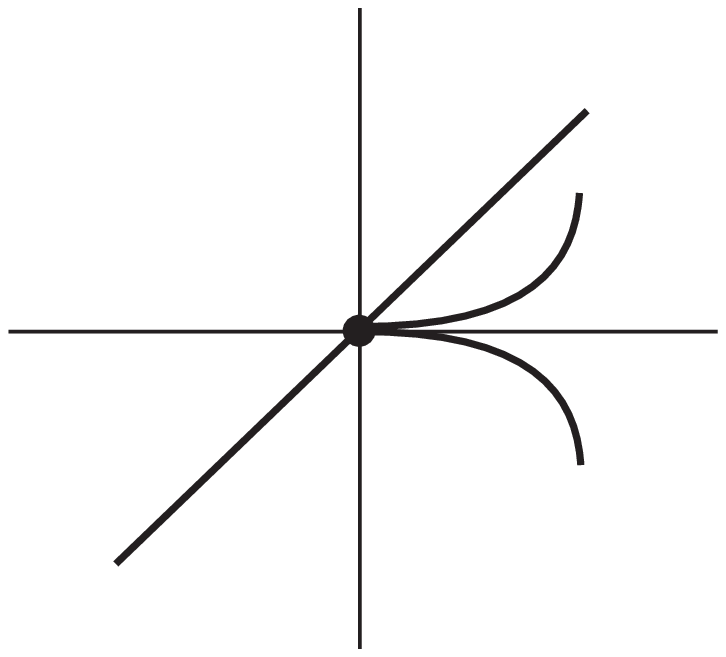}}&&\scalebox{.33}{\includegraphics{4fig3c}} \\ \\
&& Diagram Type 5: \\ &
\scalebox{.33}{\includegraphics{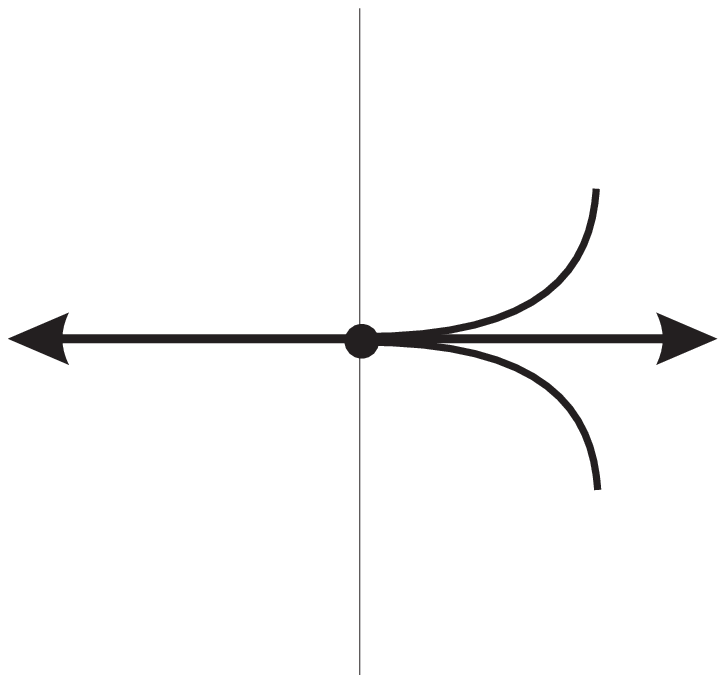}}&&\scalebox{.33}{\includegraphics{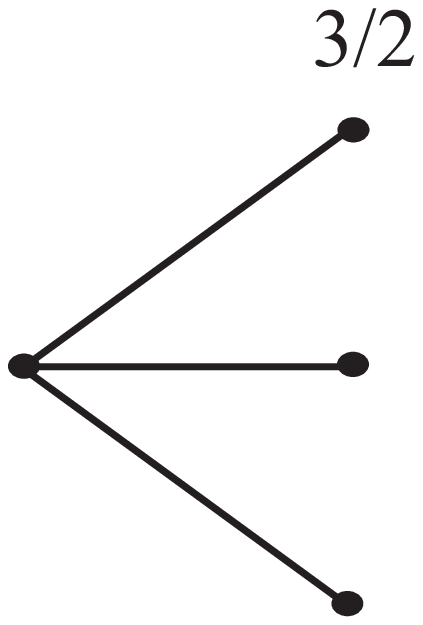}}
\\ \\\end{tabular} \end{table}
\newpage
\begin{table}[t]
\begin{tabular}{l l l l}
 && Diagram Type 6: \\ &\scalebox{.33}{\includegraphics{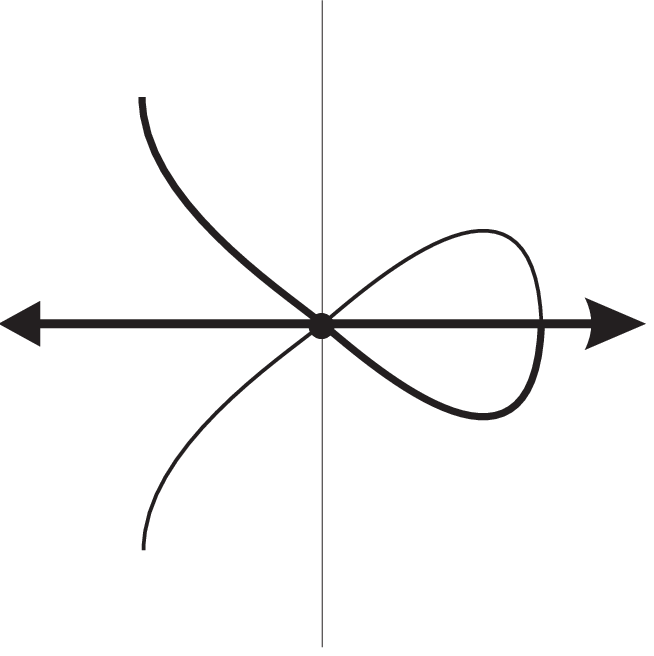}}&&\scalebox{.33}{\includegraphics{4fig9c}} \end{tabular} \\ \\
Consider the family $(y^2 - x^2 + ax^3 + bx^2y + cxy^2 + dy^3)(y - x) = 0$. \\ \\
By using Maple, we obtain\\ \\
Puiseux jets:\\
$y = x + 0x^2$\\
$y = x + \frac{1}{2} (-a-b-c-d) x^2$\\
$y = -x.$\\ \\
Condition:  $a + b + c + d \neq 0$.\\ \\
\begin{tabular}{l l l l}
 &&  Diagram Type 7: \\ &\scalebox{.33}{\includegraphics{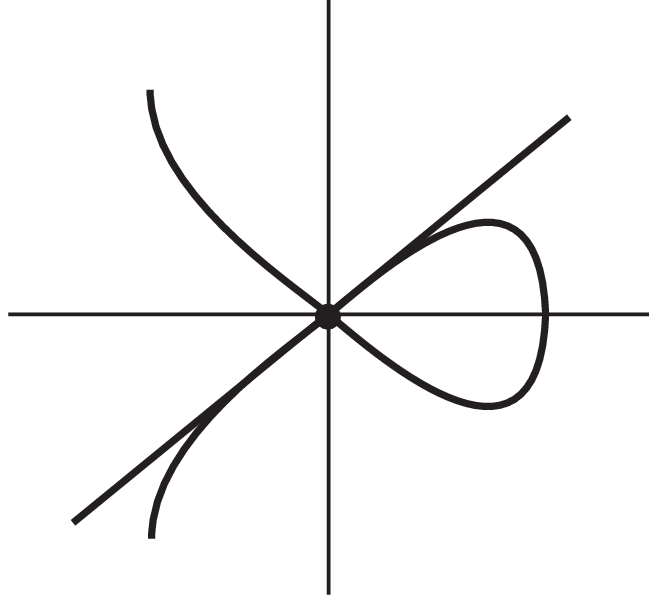}}&&\scalebox{.33}{\includegraphics{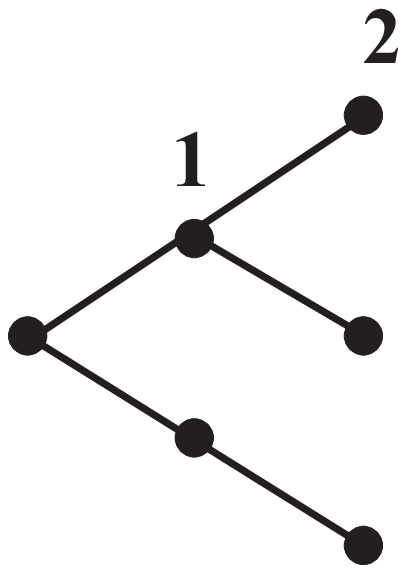}} \end{tabular} \\ \\
Case:  $a + b + c + d = 0$.\\ \\
Then the cubic is reducible, so we are done. \\ \\
\begin{tabular}{l l l l}
 && Diagram Type 8: \\ & \scalebox{.33}{\includegraphics{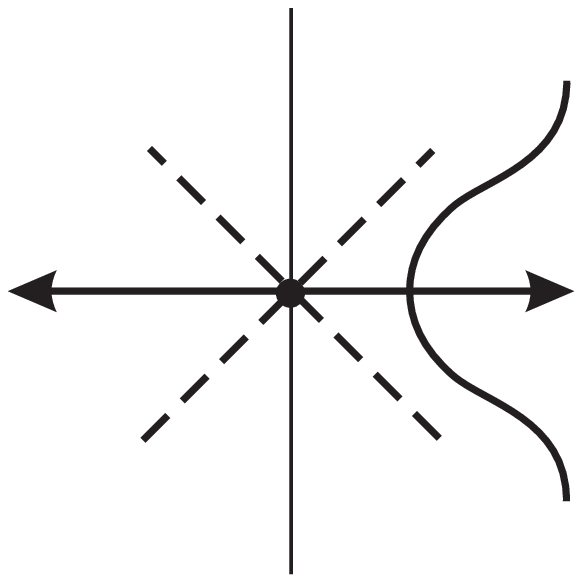}}&&\scalebox{.33}{\includegraphics{4fig11c}} \end{tabular}\\ \\
Consider the family $(y - x + ax^2 + bxy + cy^2 + dx^3 + ex^2 y + fxy^2 + gy^3)(y - x) = 0$.\\ \\
By using Maple, we obtain\\ \\
Puiseux jets:\\
$y = x + 0x^2$.\\
$y = x + (-a-b-c) x^2.$\\ \end{table} \newpage
\begin{table}[t]
Condition:  $a + b + c \neq 0$\\ \\
\begin{tabular}{l l}
 Diagram Type 9:  \\ &\scalebox{.33}{\includegraphics{4fig19c}}\end{tabular}\\ \\
Case:  $a = - b - c$.\\ \\
Puiseux jets:  \\$y = x + 0x^3$\\
$y = x + (-f-g-e-d)x^3$.\\ \\
Condition:  $f + g + e + d \neq 0$.\\ \\
\begin{tabular}{l l}
 Diagram Type 10: \\ &\scalebox{.33}{\includegraphics{4fig22c}} \end{tabular}\\ \\
Case:  $f + g + e + d = 0$.\\ \\
Then the cubic is reducible, so we are done.\\ \\
Degrees of factors: $2,2$.  \\ \\
Case:  Two distinct tangents at $(0,0)$.\\ \\
\begin{tabular}{l l l l}
 && Diagram Type 1: \\ &
\scalebox{.33}{\includegraphics{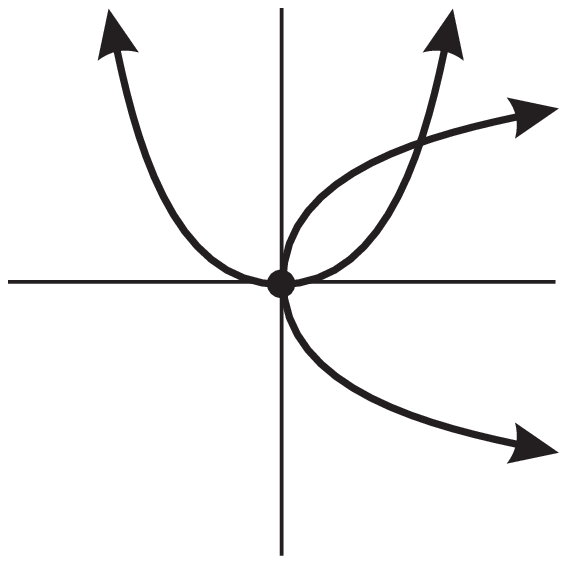}}&&\scalebox{.33}{\includegraphics{4fig2c}}
\end{tabular}
\\ \\
 $(y + ix+ ix^2) (y - ix - ix^2)= y^2 + x^2 + 2x^3 + x^4 = y^2 + x^2 (1 + x)^2$.\\ \\
\begin{tabular}{l l l l}
 && Diagram Type 2: \\ &
\scalebox{.33}{\includegraphics{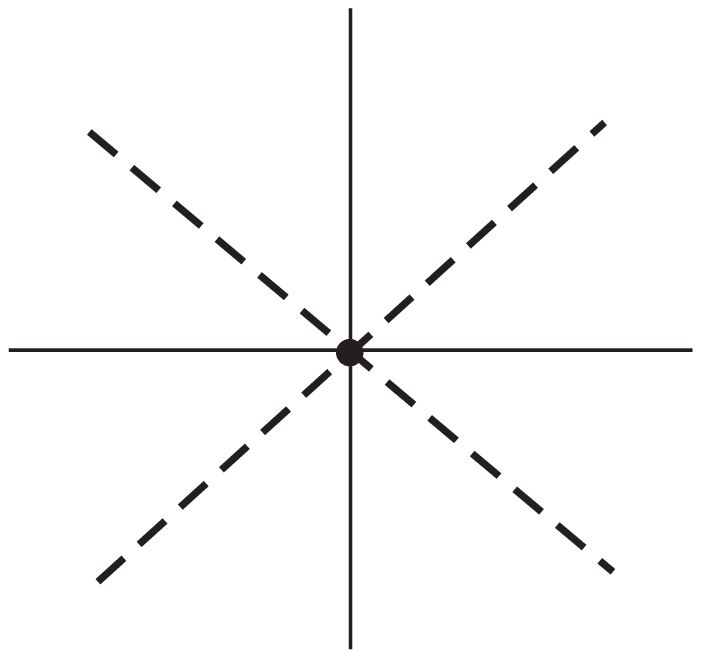}}&&\scalebox{.33}{\includegraphics{4fig15c}}
\end{tabular}\\ \\
Case:  Two tangents coincide at $(0,0)$.\\ \\ \end{table} \newpage
\begin{table}
$(y + ax^2 + bxy + cy^2)(y+dx^2 + exy + fy^2) = 0$.\\ \\
Notice that $a \neq 0$ and $d \neq 0$.
Calculation using Maple yields\\ \\
Puiseux jets:  \\ $y = -ax^2\qquad \\y = -dx^2$.\\ \\
Condition: $a \neq d$. \\ \\
\begin{tabular}{l l l l}
Diagrams: & Type 9.&& Type 11.\\
&\scalebox{.33}{\includegraphics{4fig19c}}
&&\scalebox{.33}{\includegraphics{4fig20c}}
\end{tabular}\\ \\
Case:  $a=d$.\\ \\
Puiseux jets:\\  $y = -dx^2 + edx^3 \quad \\y = - dx^2 + bdx^3$.\\ \\
Condition:  $b \neq e$.\\ \\
\begin{tabular}{l l l l}
Diagrams: & Type 10.&& Type 12.\\
&\scalebox{.33}{\includegraphics{4fig22c}}
&&\scalebox{.33}{\includegraphics{4fig23c}}
\end{tabular}\\ \\
Case:  $b = e$. \\ \\
Puiseux expansions:\\
$ y = -dx^2 +edx^3 + (-e^2d - cd^2) x^4$\\
$y = -dx^2 + edx^3 + (-e^2d - fd^2)x^4$.\\ \\
Condition: $c \neq f$.\\ \\
\begin{tabular}{l l l l}
Diagrams: & Type 13.&& Type 14.\\ &\scalebox{.33}{\includegraphics{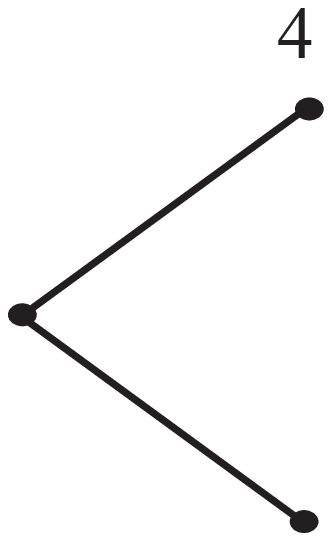}} &&\scalebox{.33}{\includegraphics{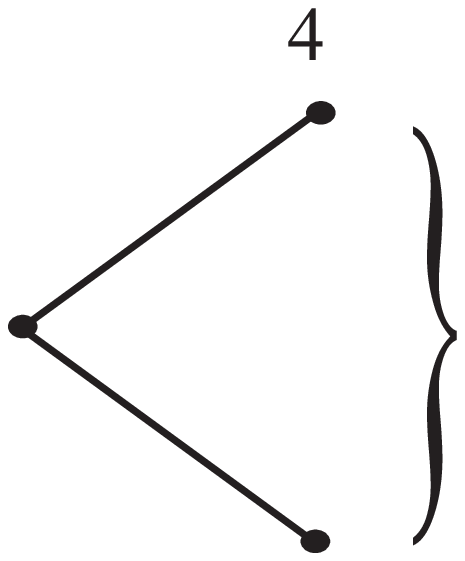}}\end{tabular}\\ \\
If $c = f$, then we just have a curve of degree $2$ with multiplicity $2$. \end{table} \newpage
\begin{table}[t]
Degrees of factors: $2,1,1$.\\ \\
\begin{tabular}{l l l l l}
  & or && Diagram Type 1: \\ \scalebox{.33}{\includegraphics{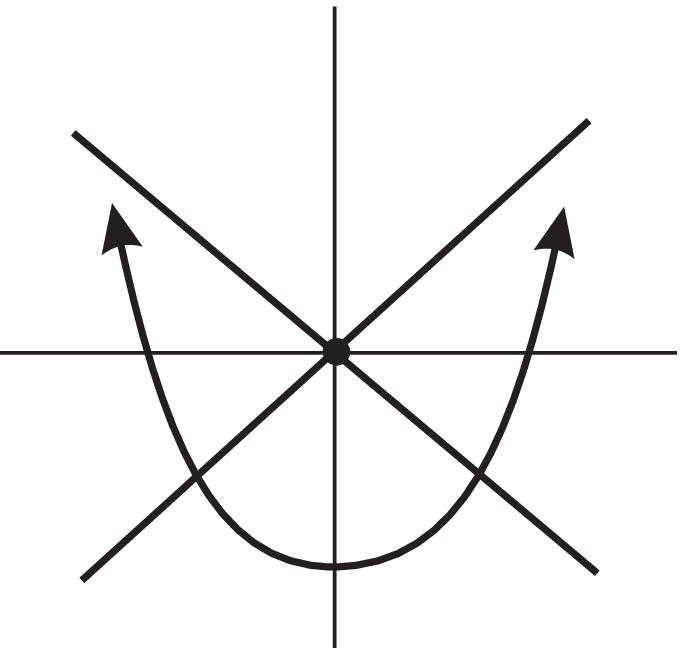}}&&\scalebox{.33}{\includegraphics{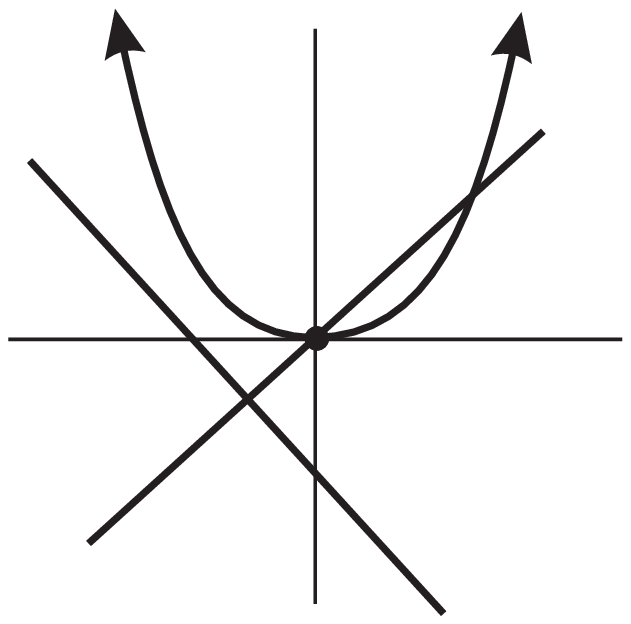}}&&\scalebox{.33}{\includegraphics{4fig2c}} \\ \\
 && Diagram Type 2: \\ \scalebox{.33}{\includegraphics{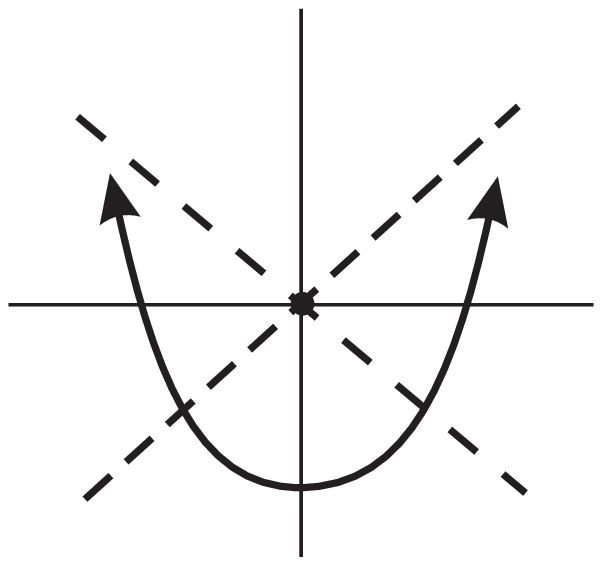}}&&&\scalebox{.33}{\includegraphics{4fig15c}} \\ \\
&& Diagram Type 9: \\ \scalebox{.33}{\includegraphics{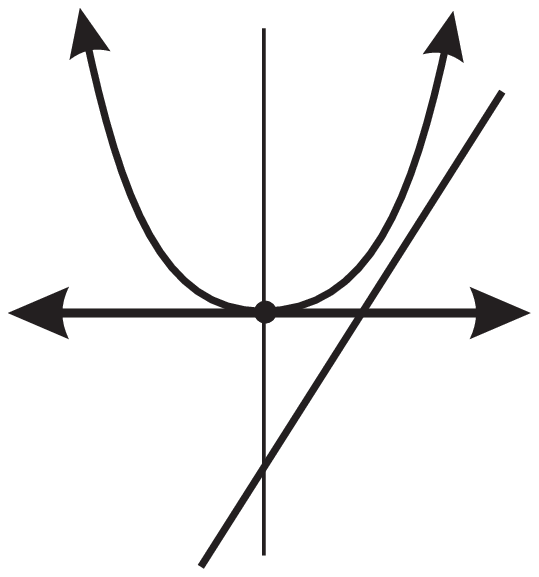}}&&&\scalebox{.33}{\includegraphics{4fig19c}} \\ \\
 && Diagram Type 6: \\ \scalebox{.33}{\includegraphics{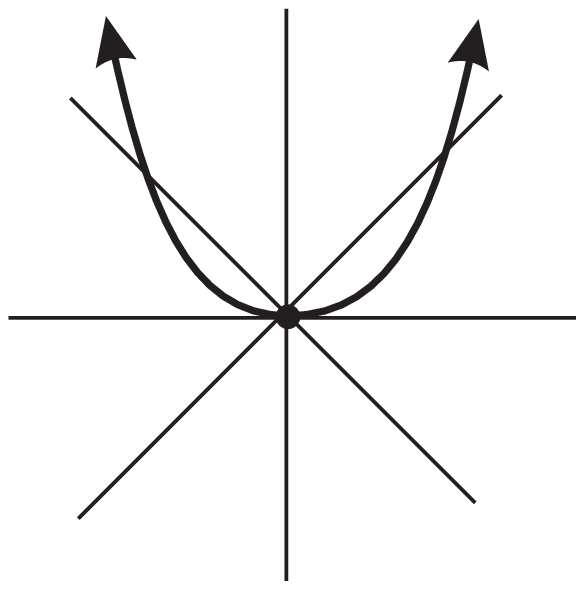}}&&&\scalebox{.33}{\includegraphics{4fig9c}} \\ \\
 && Diagram Type 8: \\ \scalebox{.33}{\includegraphics{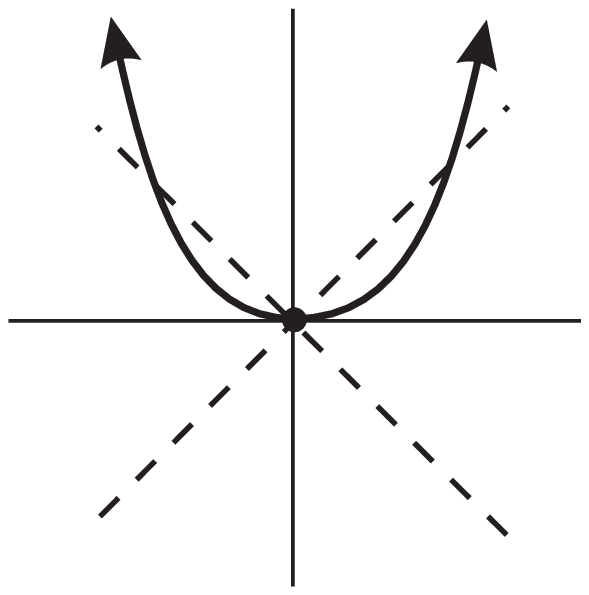}}&&&\scalebox{.33}{\includegraphics{4fig11c}} \\ \\
 && Diagram Type 7: \\ \scalebox{.33}{\includegraphics{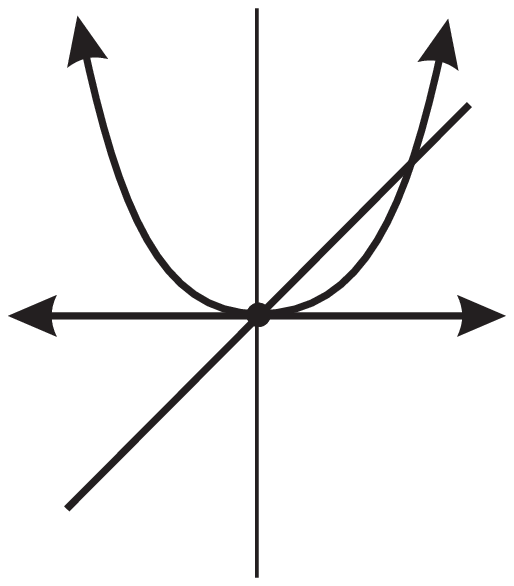}}&&&\scalebox{.33}{\includegraphics{4fig38c}} \\ \\
\end{tabular}\end{table}\newpage
\begin{table}[t]
Degrees of factors:  $1, 1, 1, 1$.\\ \\
\begin{tabular}{l l l l}
 && Diagram Type 15: \\ \scalebox{.33}{\includegraphics{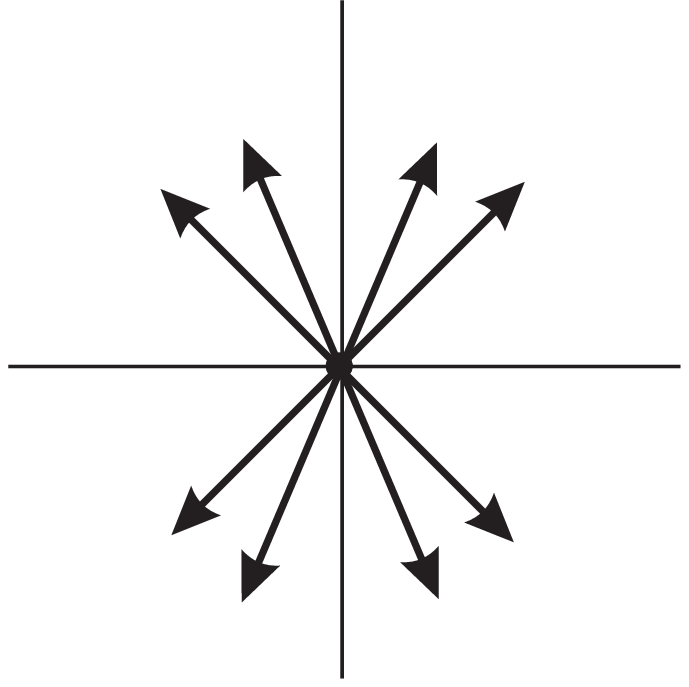}}&&&\scalebox{.33}{\includegraphics{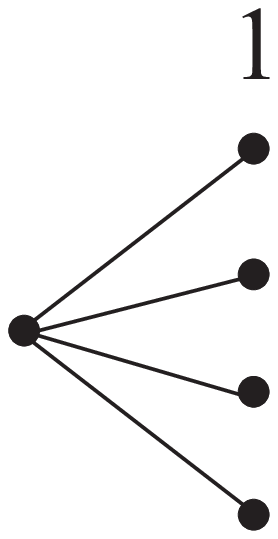}} \\ \\
 && Diagram Type 16: \\ \scalebox{.33}{\includegraphics{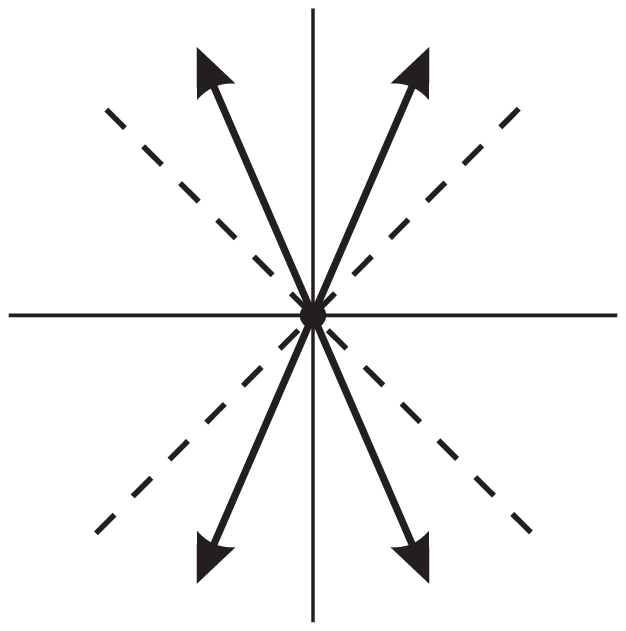}}&&&\scalebox{.33}{\includegraphics{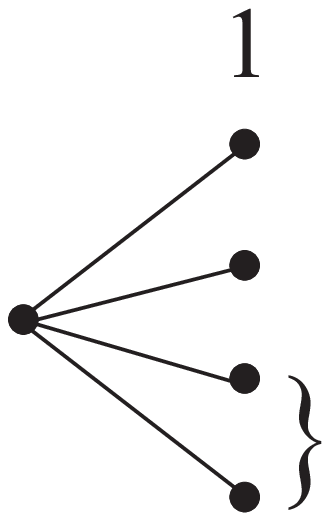}} \\ \\
 && Diagram Type 17: \\ \scalebox{.33}{\includegraphics{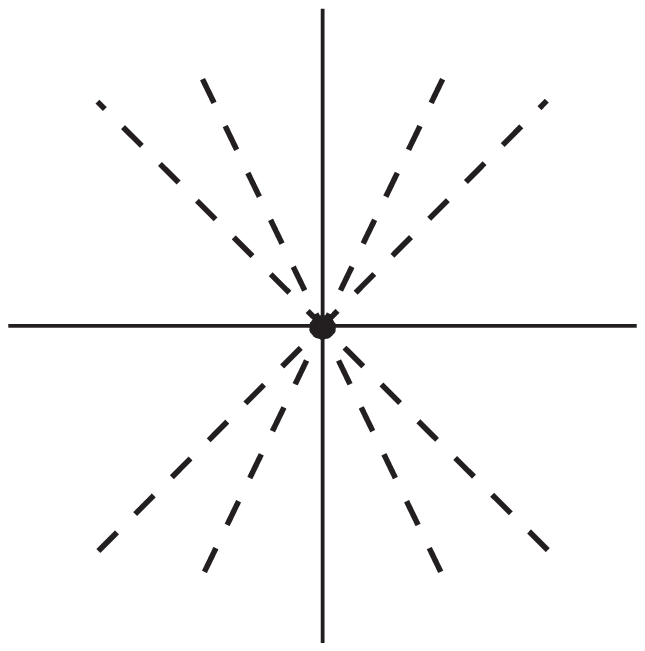}}&&&\scalebox{.33}{\includegraphics{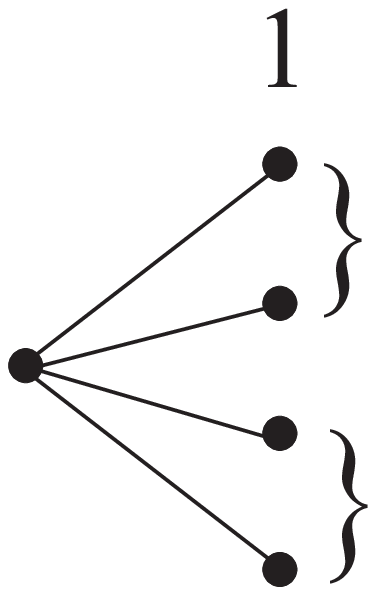}} \end{tabular} \\ \\
This completes the classification of singular point types for reducible real quartic curves. \\ \\
\section{Summary of classification.}

In this section, we summarize the classification by providing tables that show, for each singular point type, a simple example, together with a picture, the tree diagram, and the name of the singularity according to the Arnol'd notation.
(Please note that in the table of reducible curves, the example will not always perfectly match the picture.)
\end{table}  \newpage
\begin{table}[!h]
\begin{center}
\underline{Irreducible Curves} \\ \vspace{5mm}
\begin{tabular}{c|c|c|c}
Name &  Picture & Diagram  & Example\\
\hline
\begin{tabular}{c}1.  $E_6$ \\  \\  \\ \\ \end{tabular} &\scalebox{.33}{\includegraphics{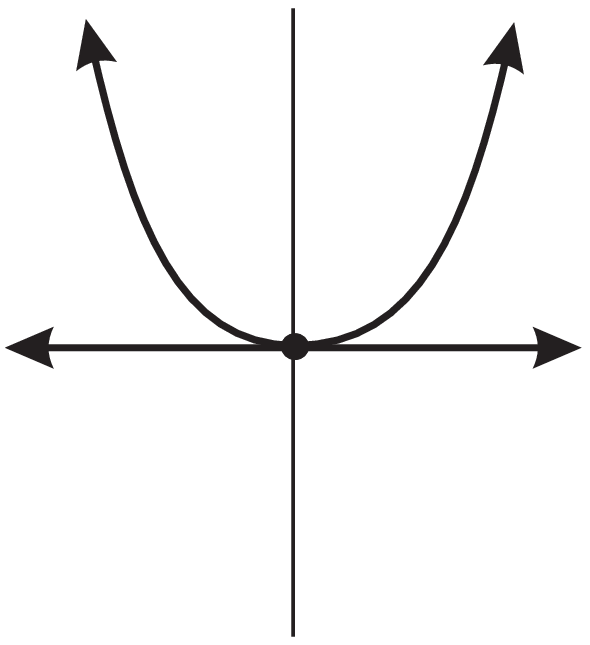}}&\scalebox{.33}{\includegraphics{4fig5c}}  &\begin{tabular}{c}
$y^3 - x^4 =0$\\ \\ \\ \\ \end{tabular} \\
\begin{tabular}{c}2.  $ D_5$ \\  \\  \\ \\ \end{tabular} &\scalebox{.33}{\includegraphics{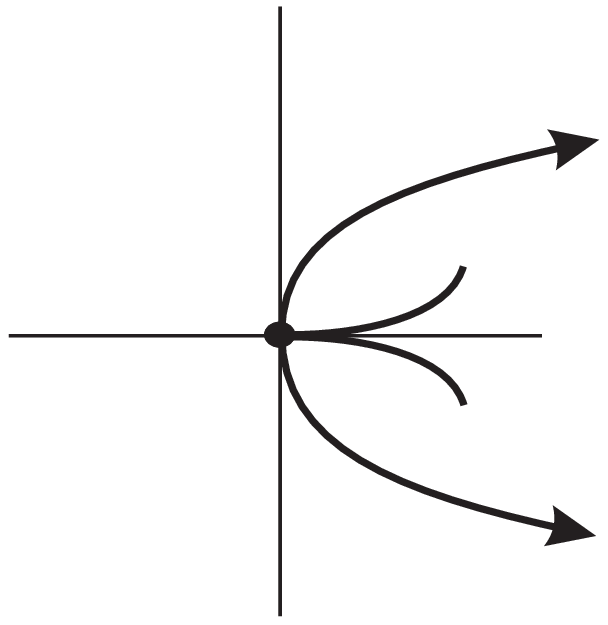}}&\scalebox{.33}{\includegraphics{4fig3c}}  &\begin{tabular}{c}$x^4+xy^2+y^3=0$\\ \\ \\ \\ \end{tabular} \\
 \begin{tabular}{c}3.  $ D_4$  \\  \\  \\ \\ \end{tabular} &\scalebox{.33}{\includegraphics{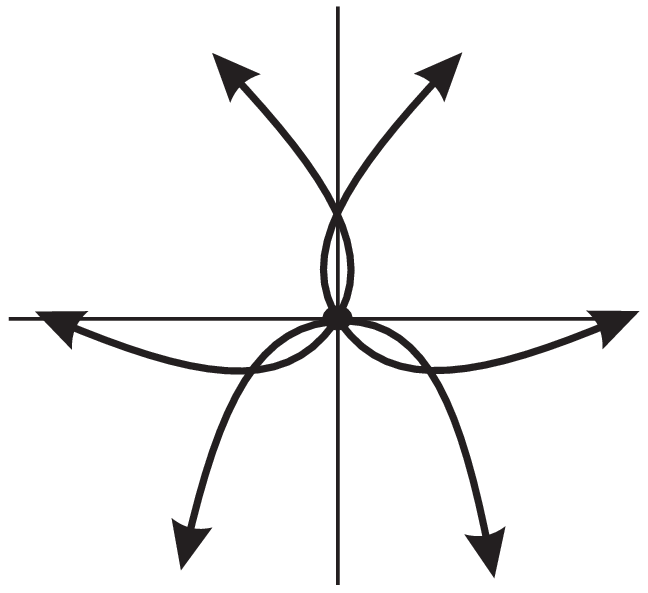}}&\scalebox{.33}{\includegraphics{4fig9c}}  &\begin{tabular}{c}
$x^4+x^2y+xy^2 +y^3=0$\\ \\ \\ \\ \end{tabular} \\
\begin{tabular}{c}4.  $ D_4^\ast$   \\  \\  \\ \\ \end{tabular} &\scalebox{.33}{\includegraphics{4fig64c}}&\scalebox{.33}{\includegraphics{4fig11c}}  &\begin{tabular}{c}
$x^4+x^2y+y^3=0$\\ \\ \\ \\ \end{tabular} \\
\begin{tabular}{c}5. $ A_1$ \\  \\  \\ \\ \end{tabular} &\scalebox{.33}{\includegraphics{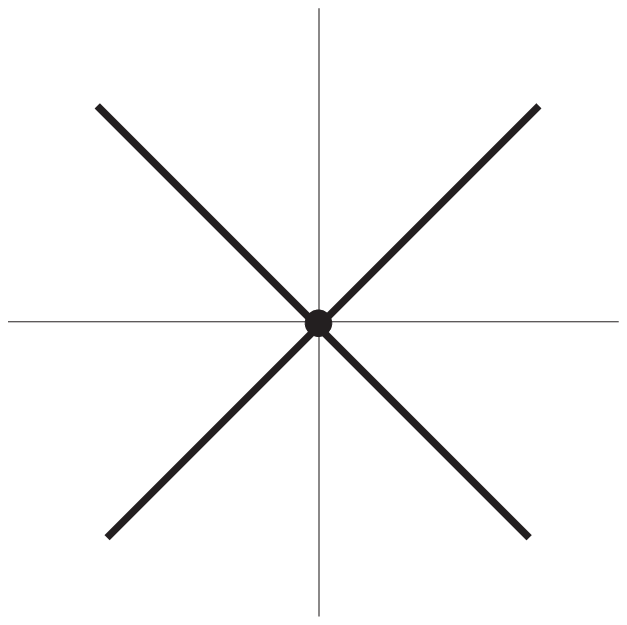}}&\scalebox{.33}{\includegraphics{4fig2c}}  &\begin{tabular}{c}
$y^2-x^2+x^4=0$\\ \\ \\ \\ \end{tabular} \\
\begin{tabular}{c}6. $ A^\ast_1$  \\  \\  \\ \\ \end{tabular} &\scalebox{.33}{\includegraphics{4fig47c}}&\scalebox{.33}{\includegraphics{4fig15c}}  &\begin{tabular}{c}
$y^2+x^2+x^4=0$\\ \\ \\ \\ \end{tabular} \\
\end{tabular}
\end{center}
\end{table}
\newpage

\begin{table}[!h]
\begin{center}
\begin{tabular}{c|c|c|c}
Name &  Picture & Diagram  & Example\\
\hline
\begin{tabular}{c}7.  $ A_2$ \\  \\  \\ \\ \end{tabular} &\scalebox{.33}{\includegraphics{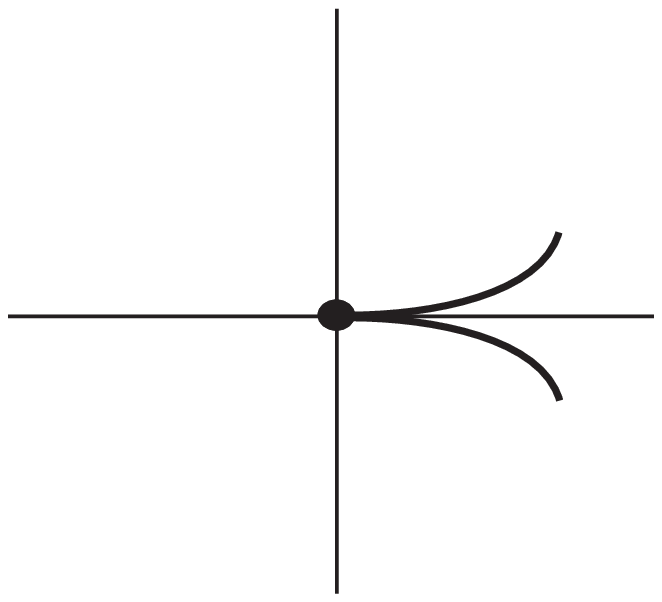}}&\scalebox{.33}{\includegraphics{4fig1c}}  &\begin{tabular}{c}
$y^2+x^3+x^4=0$\\ \\ \\ \\ \end{tabular} \\
\begin{tabular}{c}8.  $ A_3$ \\  \\  \\ \\ \end{tabular} &\scalebox{.33}{\includegraphics{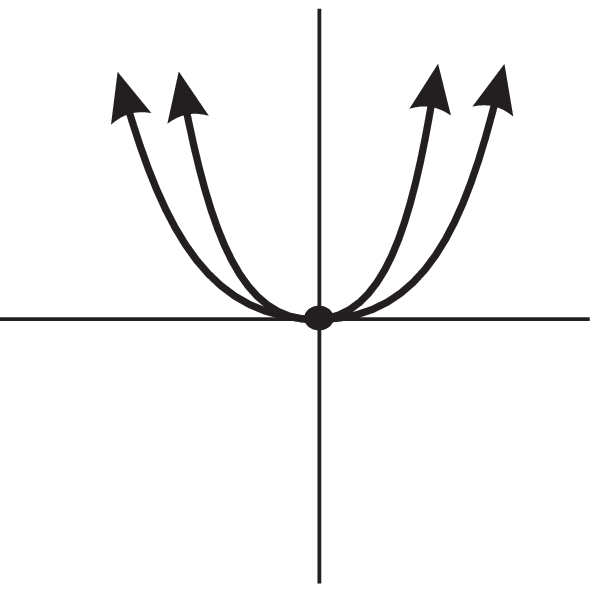}}&\scalebox{.33}{\includegraphics{4fig19c}}  &\begin{tabular}{c}$y^2-x^4+y^3=0$\\ \\ \\ \\ \end{tabular} \\
\begin{tabular}{c}9.  $ A^\ast_3$ \\  \\  \\ \\ \end{tabular} &\scalebox{.33}{\includegraphics{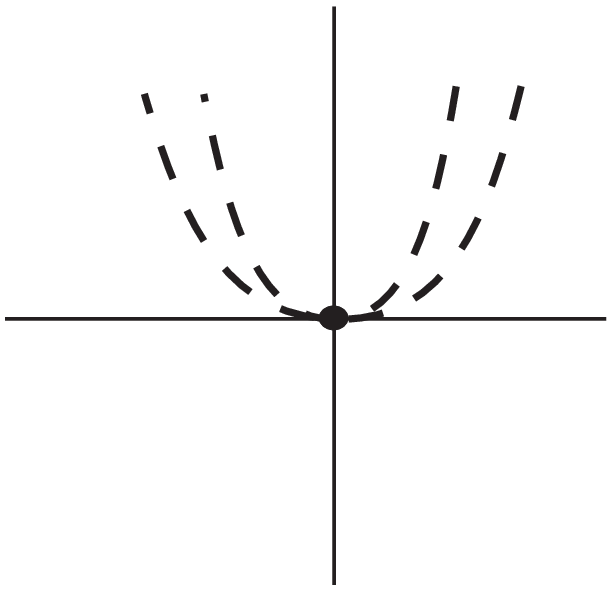}}&\scalebox{.33}{\includegraphics{4fig20c}}  &\begin{tabular}{c}$y^2+x^4+y^3=0$\\ \\ \\ \\ \end{tabular} \\
\begin{tabular}{c}10.  $ A_4$ \\  \\  \\ \\ \end{tabular} &\scalebox{.33}{\includegraphics{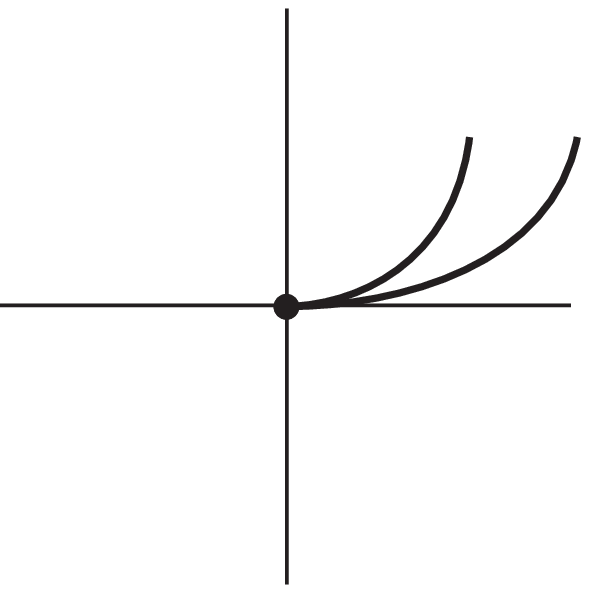}}&\scalebox{.33}{\includegraphics{4fig21c}}  &\begin{tabular}{c}$y^2+2x^2y+x^4+x^3y=0$\\ \\ \\ \\ \end{tabular} \\
\begin{tabular}{c}11.  $ A_5$ \\  \\  \\ \\ \end{tabular} &\scalebox{.33}{\includegraphics{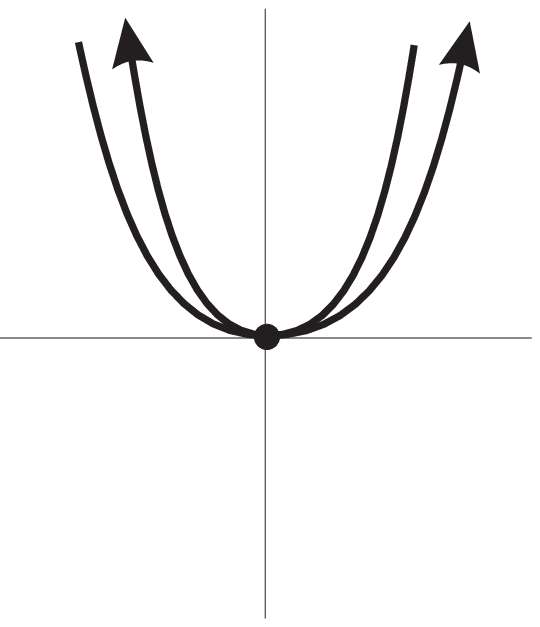}}&\scalebox{.33}{\includegraphics{4fig22c}}  &\begin{tabular}{c}$y^2+2x^2y+x^4+y^3=0$\\ \\ \\ \\ \end{tabular} \\
\begin{tabular}{c}12.  $ A^\ast_5$ \\  \\  \\ \\ \end{tabular} &\scalebox{.33}{\includegraphics{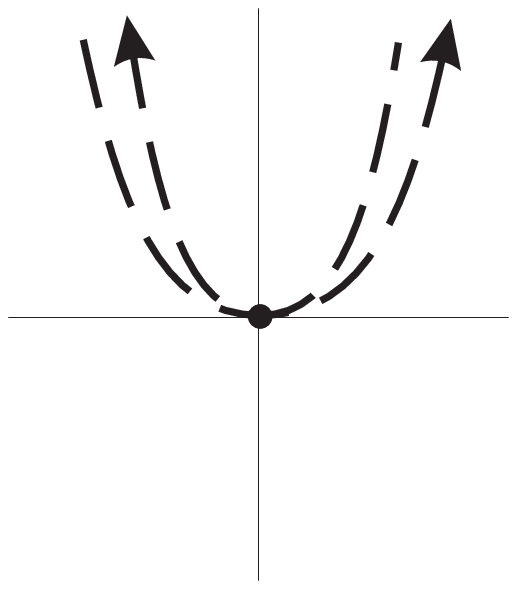}}&\scalebox{.33}{\includegraphics{4fig23c}}  &\begin{tabular}{c}$y^2+2x^2y+x^4-y^3=0$\\ \\ \\ \\ \end{tabular} \\

\end{tabular}
\end{center}
\end{table}

\newpage
\begin{table}[!h]
\begin{center}
\begin{tabular}{c|c|c|c}
Name &  Picture & Diagram  & Example\\
\hline

\begin{tabular}{c}13.  $ A_6$ \\  \\  \\ \\ \end{tabular} &\scalebox{.33}{\includegraphics{4fig80c}}&\scalebox{.33}{\includegraphics{4fig24c}}  &\begin{tabular}{c}$y^2+2x^2y+x^4+x^2y^2+\frac{1}{4}y^4+y^3=0$\\ \\ \\ \\ \end{tabular} \\
\end{tabular}
\vspace{3mm} \\
\underline{Reducible Curves}\\
\vspace{3mm}
\begin{tabular}{c|c|c|c}
Name &  Picture & Diagram  & Example\\
\hline
\begin{tabular}{c}1.  $A_1$ \\  \\  \\ \\ \end{tabular} &\scalebox{.33}{\includegraphics{4figA1c}}&\scalebox{.33}{\includegraphics{4fig2c}}  &\begin{tabular}{c}$(y-1)(y-2)(y-x)(y+x)=0$\\ \\ \\ \\ \end{tabular} \\
\begin{tabular}{c}2.  $A_1^\ast$ \\  \\  \\ \\ \end{tabular} &\scalebox{.33}{\includegraphics{4fig47c}}&\scalebox{.33}{\includegraphics{4fig15c}}  &\begin{tabular}{c}$(y-1)(y-2)(x^2+y^2)=0$\\ \\ \\ \\ \end{tabular} \\
\begin{tabular}{c}3.  $A_2$ \\  \\  \\ \\ \end{tabular} &\scalebox{.33}{\includegraphics{4fig77c}}&\scalebox{.33}{\includegraphics{4fig1c}}  &\begin{tabular}{c}$(y-1)(y^2-x^3)=0$\\ \\ \\ \\ \end{tabular} \\
\begin{tabular}{c}4.  $D_5$ \\  \\  \\ \\ \end{tabular} &\scalebox{.33}{\includegraphics{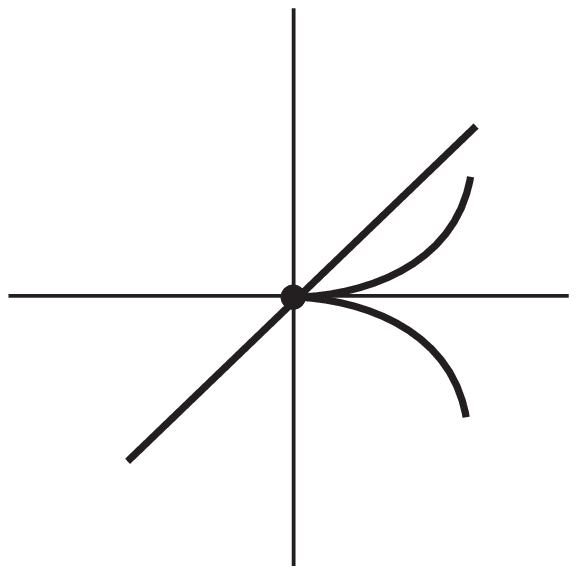}}&\scalebox{.33}{\includegraphics{4fig3c}}  &\begin{tabular}{c}$(y-x)(y^2-x^3)=0$\\ \\ \\ \\ \end{tabular} \\
\begin{tabular}{c}5.  $E_7$ \\  \\  \\ \\ \end{tabular} &\scalebox{.33}{\includegraphics{4fig33c}}&\scalebox{.33}{\includegraphics{4figmissing_c}}  &\begin{tabular}{c}$y(y^2-x^3)=0$\\ \\ \\ \\ \end{tabular} \\

\end{tabular}
\end{center}
\end{table}

\newpage
\begin{table}[!h]
\begin{center}
\begin{tabular}{c|c|c|c}
Name &  Picture & Diagram  & Example\\
\hline
\begin{tabular}{c}6.  $D_4$ \\  \\  \\ \\ \end{tabular} &\scalebox{.33}{\includegraphics{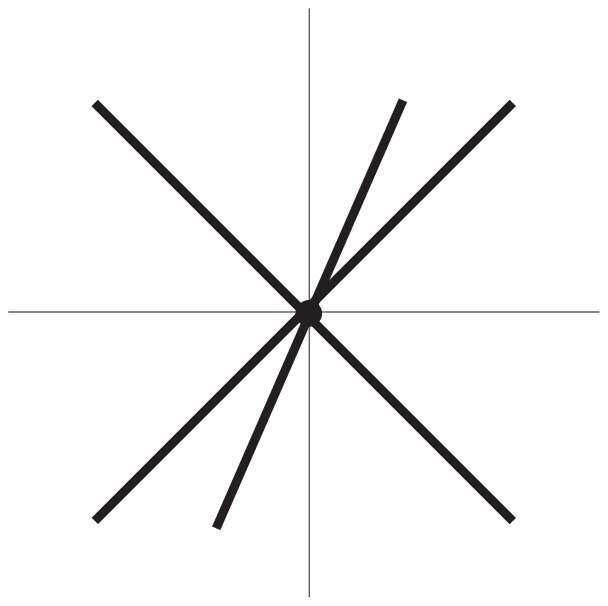}}&\scalebox{.33}{\includegraphics{4fig9c}}  &\begin{tabular}{c}$(y-x)(y+x)(y-x^2)=0$\\ \\ \\ \\ \end{tabular} \\
\begin{tabular}{c}7.  $D_6$ \\  \\  \\ \\ \end{tabular} &\scalebox{.33}{\includegraphics{4fig66c}}&\scalebox{.33}{\includegraphics{4fig38c}}  &\begin{tabular}{c}$y(y-x)(y-x^2)=0$\\ \\ \\ \\ \end{tabular} \\
\begin{tabular}{c}8.  $D^\ast_4$ \\  \\  \\ \\ \end{tabular} &\scalebox{.33}{\includegraphics{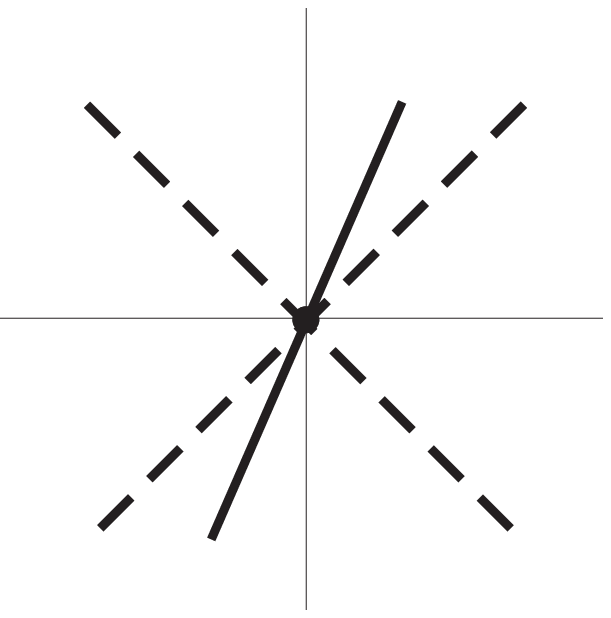}}&\scalebox{.33}{\includegraphics{4fig11c}}  &\begin{tabular}{c}$(y-x)(y^2+x^2-x^3)=0$\\ \\ \\ \\ \end{tabular} \\
\begin{tabular}{c}9.  $A_3$ \\  \\  \\ \\ \end{tabular} &\scalebox{.33}{\includegraphics{4fig78c}}&\scalebox{.33}{\includegraphics{4fig19c}}  &\begin{tabular}{c}$(y-x^2)(y+x^2)=0$\\ \\ \\ \\ \end{tabular} \\
\begin{tabular}{c}10.  $A_3^\ast$ \\  \\  \\ \\ \end{tabular} &\scalebox{.33}{\includegraphics{4fig79c}}&\scalebox{.33}{\includegraphics{4fig20c}}  &\begin{tabular}{c}$y^2+x^4=0$\\ (Reducible over $\mathbb{C}$) \\ \\ \\ \end{tabular} \\
\begin{tabular}{c}11.  $A_5$ \\  \\  \\ \\ \end{tabular} &\scalebox{.33}{\includegraphics{4figA5c}}&\scalebox{.33}{\includegraphics{4fig22c}}  &\begin{tabular}{c}$(y+x^2)(y+x^2+xy)=0$\\ \\ \\ \\ \end{tabular} \\
\end{tabular}
\end{center}
\end{table}
\newpage
\begin{table}[!h]
\begin{center}
\begin{tabular}{c|c|c|c}
Name &  Picture & Diagram  & Example\\
\hline
\begin{tabular}{c}12.  $A_5^\ast$ \\  \\  \\ \\ \end{tabular} &\scalebox{.33}{\includegraphics{4figA5xc}}&\scalebox{.33}{\includegraphics{4fig23c}}  &\begin{tabular}{c}$x^4+2x^2y+y^2x^2+y^2=0$\\(Reducible over $\mathbb{C}$) \\ \\ \\ \end{tabular} \\
\begin{tabular}{c}13.  $A_7$ \\  \\  \\ \\ \end{tabular} &\scalebox{.33}{\includegraphics{4fig78c}}&\scalebox{.33}{\includegraphics{4fig53c}}  &\begin{tabular}{c}$(y+x^2+xy)(y+x^2+xy+y^2)=0$\\ \\ \\ \\ \end{tabular} \\
\begin{tabular}{c}14.  $A_7^\ast$ \\  \\  \\ \\ \end{tabular} &\scalebox{.33}{\includegraphics{4fig79c}}&\scalebox{.33}{\includegraphics{4fig54c}}  &\begin{tabular}{c}$x^4+2x^3y+2x^2y+y^2x^2+2xy^2+y^4+y^2=0$\\ (Reducible over $\mathbb{C}$)\\ \\ \\ \end{tabular} \\
\begin{tabular}{c}15.  $X_9$ \\  \\  \\ \\ \end{tabular} &\scalebox{.33}{\includegraphics{4fig68c}}&\scalebox{.33}{\includegraphics{4fig69c}}  &\begin{tabular}{c}$(y-x)(y+x)(y-2x)(y+2x)=0$\\ \\ \\ \\ \end{tabular} \\
\begin{tabular}{c}16.  $X^\ast_9$ \\  \\  \\ \\ \end{tabular} &\scalebox{.33}{\includegraphics{4fig70c}}&\scalebox{.33}{\includegraphics{4fig71c}}  &\begin{tabular}{c}$(y-x)(y+x)(x^2+y^2)=0$\\ \\ \\ \\ \end{tabular} \\
\begin{tabular}{c}17.  $X^{\ast\ast}_9$ \\  \\  \\ \\ \end{tabular} &\scalebox{.33}{\includegraphics{4fig72c}}&\scalebox{.33}{\includegraphics{4fig73c}}  &\begin{tabular}{c}$(x^2+y^2)(x^2+4y^2)=0$\\ \\ \\ \\ \end{tabular} \\
\end{tabular}
\end{center}
\end{table}
\newpage

\begin{table}[!t]

\end{table}
\begin{table}[h]
David Weinberg: Department of Mathematics and Statistics, Texas
Tech University, Lubbock, TX 79409-1042

Nicholas Willis: Mathematics and Computer Science Department,
Whitworth College, Spokane, WA 99251
\end{table}

\end{document}